\documentclass{article}

\usepackage{amsmath}
\usepackage{amsthm}
\usepackage[psamsfonts]{amsfonts}
\usepackage{amssymb}
\usepackage[all]{xy}

\DeclareMathSymbol{\rightrightarrows}  {\mathrel}{AMSa}{"13}

\def\Ob{\operatorname{Ob}}

\def\Mor{\operatorname{Mor}}
\def\St{\operatorname{St}}

\def\Ho{\operatorname{Ho}}

\def\Pre{\operatorname{Pre}}

\def\Ext{\operatorname{Ext}}

\def\pb{\operatorname{pb}}

\catcode`\@=11
\def\varholim@#1#2{\mathop{\vtop{\ialign{##\crcr
 \hfil$#1\m@th\operator@font holim$\hfil\crcr
 \noalign{\nointerlineskip\kern\ex@}#2#1\crcr
 \noalign{\nointerlineskip\kern-\ex@}\crcr}}}}
\def\hocolim{\mathpalette\varholim@\rightarrowfill@} 
\def\hoinvlim{\mathpalette\varholim@\leftarrowfill@}
\catcode`\@=\active

\newtheorem{theorem}{Theorem}
\newtheorem{lemma}[theorem]{Lemma}
\newtheorem{corollary}[theorem]{Corollary}

\theoremstyle{definition}

\newtheorem{example}[theorem]{Example}

\newtheorem{remark}[theorem]{Remark}

\begin{document}

\title{Cocycle categories}
                                                                                \author{J.F. Jardine\thanks{This research was supported by NSERC.}}
 
 
\date{March 29, 2006}

\maketitle

\section*{Introduction}

Suppose that $G$ is a sheaf of groups on a space $X$ and that
$U_{\alpha} \subset X$ is an open covering. Then a cocycle for the
covering is traditionally defined to be a family of elements
$g_{\alpha\beta} \in G(U_{\alpha} \cap U_{\beta})$ such that
$g_{\alpha\alpha} = e$ and $g_{\alpha\beta}g_{\beta\gamma} =
g_{\alpha\gamma}$ when all elements are restricted to the group
$G(U_{\alpha}\cap U_{\beta} \cap U_{\gamma})$. 

A more compact way of
saying this is to assert that such a cocycle is a map of simplicial
sheaves $C(U) \to BG$ on the space $X$, where $C(U)$ is the \v{C}ech
resolution associated to the covering family $\{ U_{\alpha} \}$. The
canonical map $C(U) \to \ast$ is a local weak equivalence of
simplicial sheaves, and is a fibration in each section since $C(U)$ is
actually the nerve of a groupoid --- the map $C(U) \to \ast$ is
therefore a hypercover, which is most properly defined to be a map of
simplicial sheaves (or presheaves) which is a Kan fibration and a weak
equivalence in each stalk. Every cocycle in the traditional sense
therefore determines a picture of simplicial sheaf morphisms
\begin{equation*}
\ast \leftarrow C(U) \to BG.
\end{equation*}
where the canonical map $C(U) \to \ast$ is a hypercover.

More generally, it has been known for a long time \cite{Br}, \cite{J1}
that the locally fibrant simplicial sheaves (ie. simplicial sheaves
which are Kan complexes in each stalk) have a partial homotopy theory
which is good enough to give a calculus of fractions approach to
formally inverting local weak equivalences, and thereby giving a
construction of the homotopy category for simplicial
sheaves. Specifically, one defines morphisms $[X,Y]$ in the homotopy
category by setting
\begin{equation}\label{eq 00}
[X,Y] = \varinjlim_{[\pi]:V \to X} \pi(V,Y),
\end{equation}
where the filtered colimit is indexed over simplicial homotopy classes
of hypercovers $\pi: V \to X$, and $\pi(V,Y)$ denotes simplicial
homotopy classes of maps $V \to Y$. Thus, morphisms in the homotopy
category are represented by pictures
\begin{equation}\label{eq 0}
X \xleftarrow{\pi} V \to Y,
\end{equation}
where $\pi$ is a hypercover.
The relation (\ref{eq 00}) is usually called the Generalized Verdier
Hypercovering Theorem, and it is the historical starting point for the
homotopy theory of simplicial sheaves.  

Much has transpired in the intervening years. We now know that there
is a plethora of Quillen model structures for simplicial sheaves and
simplicial presheaves on all small Grothendieck sites \cite{J6.25}, all of which determine the same homotopy category. In
the first of these structures \cite{Joy}, \cite{J2}, the monomorphisms
are the cofibrations and the weak equivalences are the local weak
equivalences (ie. stalkwise weak equivalences in the presence of
enough stalks), and then the homotopy categories for simplicial
sheaves and presheaves (which are equivalent) are constructed by
methods introduced by Quillen. The morphisms $[X,Y]$ in the homotopy
category of simplicial sheaves coincide with those specified by
(\ref{eq 00}) if $X$ and $Y$ are locally fibrant.  These basic model
structures also inherit many of the good properties from simplicial
sets, including right properness, which means that weak equivalences
are stable under pullback along fibrations. Weak equivalences of
simplicial presheaves and sheaves are also closed under finite
products, essentially because the same statement holds for simplicial
sets.

One still needs a cocycle or hypercover theory, because it is involved
in the proofs of many of the standard homotopy classification theorems
for sheaf cohomology, such as the identification of sheaf cohomology
with morphisms $[\ast,K(A,n)]$ in the homotopy category.  It turns
out, however, that the traditional requirement that the map $\pi$ in
(\ref{eq 0}) should be something like a fibration in formal
manipulations based on the generalized Verdier hypercover theorem is
usually quite awkward in practice. This has been especially apparent
in all attempts to convert $n$-types to more algebraic objects.

The basic point of this paper is that there is a better approach to
defining and manipulating cocycles, which essentially starts with
removing the fibration condition on weak equivalence $\pi$ in (\ref{eq
0}). 

A cocycle from $X$ to $Y$ is defined in this paper to be a pair of maps
\begin{equation}\label{eq 01}
X \xleftarrow{f} Z \xrightarrow{g} Y
\end{equation}
where $f$ is a weak equivalence. A morphism of cocycles is the obvious
thing, namely a commutative diagram
\begin{equation*}
\xymatrix@R=5pt{
& Z \ar[dl]_{f} \ar[dr]^{g} \ar[dd] & \\
X && Y \\
& Z' \ar[ul]^{f'} \ar[ur]_{g'} &
}
\end{equation*}
The cocycle category is denoted by $H(X,Y)$. Then the basic point is
this: there is an obvious function
\begin{equation*}
\phi: \pi_{0}H(X,Y) \to [X,Y]
\end{equation*}
defined by sending the cocycle (\ref{eq 01}) to the morphism $gf^{-1}$
in the homotopy category. Then Theorem \ref{th 2} of this paper
implies that this function $\phi$ is a bijection.  

Cocycle categories and the function $\phi$ are defined in great
generality.  The fact that $\phi$ is a bijection (Theorem \ref{th 2})
holds in the rather common setting of a right proper closed model
structure in which weak equivalences are preserved by finite
products. The formal definition of cocycle categories and their basic
properties appear in the first section of this paper. The overall
theory is rather easy to demonstrate.

The subsequent sections of this paper are taken up with a tour of the
applications of this theory.  Some of these applications are well
known, and are given here with simple new proofs.  This general
approach to cocycles is also implicated in several recent results in
non-abelian cohomology theory, which are displayed here.

Section 2 gives a quick (and hypercover-free) demonstration of the
homotopy classification theorem for $G$-torsors for a sheaf of groups
$G$ --- this is Theorem \ref{th 5}. It is also shown that the ideas
behind the proof of Theorem \ref{th 5} admit great generalization:
examples include the definition and homotopy classification of torsors
for a presheaf of categories enriched in simplicial sets (Theorem
\ref{th 6}), the explicit construction of the stack completion of a
sheaf of groupoids given in Section 2.3, and the calculation of the
morphism set $[\ast,\hocolim_{G} X]$ for a diagram $X$ of simplicial
presheaves defined on a presheaf of groupoids $G$ in Theorem 
\ref{th 7.5}. 

Section 3 gives a new demonstration of homotopy classification
theorem for abelian sheaf cohomology (Corollary \ref{cor 9}). Section
4 gives a new description of the homotopy classification theorem for
group extensions in terms of cocycles taking values in $2$-groupoids
(Theorem \ref{th 10}), and in Section 5 we discuss but do not prove a
classification theorem for gerbes up to local equivalence as path
components of a suitable cocycle category. This last result is Theorem
\ref{th 11} --- it is proved in \cite{J9}. The cocycles appearing in
the proofs of Theorems \ref{th 5}, \ref{th 10} and \ref{th 11} are
canonically defined.

\section{Cocycles}

Suppose that $\mathcal{M}$ is a right proper closed model category,
and suppose further that that the weak equivalences of $\mathcal{M}$
are closed under finite products.

The assertion that $\mathcal{M}$ is a closed model category means
that $\mathcal{M}$ has all finite limits and colimits ({\bf
CM1}), and that the class of morphisms of $\mathcal{M}$ contain three
subclasses, namely weak equivalences, fibrations and cofibrations
which satisfy some properties. These properties include the two of
three condition for weak equivalences ({\bf CM2}: if any two of the
maps $f$, $g$ or $g\cdot f$ are weak equivalences, then so is the
third), the requirement that all three classes of maps are closed
under retraction ({\bf CM3}), and the factorization axiom ({\bf CM5})
which asserts that any map $f$ has factorizations $f=qj=pi$ where $q$
is a trivial fibration and $j$ is a cofibration, and $p$ is a
fibration and $i$ is a trivial cofibration. Note, for example, that a
trivial fibration is a fibration and a weak equivalence ---
``trivial'' things are always weak equivalences. Finally,
$\mathcal{M}$ should satisfy the lifting axiom ({\bf CM4}) which says
that in any solid arrow diagram
\begin{equation*}
\xymatrix{
A \ar[r] \ar[d]_{i} & X \ar[d]^{p} \\
B \ar[r] \ar@{.>}[ur] & Y
}
\end{equation*}
where $p$ is a fibration and $i$ is a cofibration, the dotted
arrow exists making the diagram commute if either $i$ or $p$ is
trivial.

A model category $\mathcal{M}$ is said to be right proper if weak
equivalences are closed under pullback along fibrations. 

Not every model category is right proper, but right proper model
structures are fairly common: examples include topological spaces,
simplicial sets, spectra, simplicial presheaves, presheaves of
spectra, and certain ``good'' localizations such as the motivic and
motivic stable model categories.  In all of these examples as well,
weak equivalences are closed under finite products, meaning that if
$f: X \to Y$ is a weak equivalence and $Z$ is any other object, then
the map $f \times 1: X \times Z \to Y \times Z$ is a weak equivalence.
\medskip

Suppose that $X$ and $Y$ are objects of $\mathcal{M}$. Let
$H(X,Y)$ be the 
category whose objects are all pairs of maps $(f,g)$
\begin{equation*}
X \xleftarrow{f} Z \xrightarrow{g} Y
\end{equation*}
such that $f$ is a weak equivalence. A morphism $\alpha: (f,g)
\to (f',g')$ of $H(X,Y)$ is a commutative diagram
\begin{equation*}
\xymatrix@R=5pt{
& Z \ar[dl]_{f} \ar[dr]^{g} \ar[dd]^{\alpha} & \\
X && Y \\
& Z' \ar[ul]^{f'} \ar[ur]_{g'} &
}
\end{equation*}
$H(X,Y)$ is the {\it category of cocycles}
from $X$ to $Y$.

\begin{example}
Suppose that $V \to \ast$ is a sheaf epimorphism (possibly arising
from a covering) and that $G$ is a sheaf of groups. Traditional
cocycles for the underlying site with coefficients in $G$ can be
interpreted as simplicial sheaf maps
\begin{equation*}
\ast \leftarrow C(V) \to BG,
\end{equation*}
where $C(V)$ is the \v{C}ech resolution for the cover, and the
canonical map $C(V) \to \ast$ is a local weak equivalence of
simplicial sheaves.  All cocycles of this type therefore represent
objects of the cocycle category $H(\ast,BG)$ in simplicial
sheaves. This is a motivating example for the present definition of a
cocycle.

More generally, it has been known for some time \cite{J1} that
morphisms $X \to Y$ in the homotopy category of simplicial sheaves
can be represented by pairs of morphisms
\begin{equation*} 
X \xleftarrow{\pi} U \to Y
\end{equation*}
where $\pi$ is a hypercover, or locally trivial fibration. This is
provided that $X$ and $Y$ are locally fibrant. Such pictures are
members of the cocycle category $H(X,Y)$ in simplicial sheaves.
\end{example}

Generally, write $\pi_{0}H(X,Y)$ for the class of path components of
$H(X,Y)$, and write $[X,Y]$ for the set of morphisms from $X$ to $Y$
in the homotopy category $\Ho(\mathcal{M})$ of $\mathcal{M}$. Recall
that $\Ho(\mathcal{M})$ is constructed from $\mathcal{M}$ by formally
inverting the weak equivalences. Then one sees that there is a
function
\begin{equation*}
\phi: \pi_{0}H(X,Y) \to [X,Y]
\end{equation*}
which is defined by $(f,g) \mapsto g\cdot f^{-1}$.

\begin{theorem}\label{th 2}
Suppose that $\mathcal{M}$ is a right proper closed model category for
which the class of weak equivalences is closed under finite products.
Then the function $\phi: \pi_{0}H(X,Y) \to [X,Y]$ is a bijection for all
$X$ and $Y$.
\end{theorem}

Suppose that the maps $f,g: X \to
Y$ are left homotopic. Then there is a commutative diagram
\begin{equation*}
\xymatrix{
& X \ar[dr]^{f} \ar[dl]_{1} \ar[d]_{d_{0}} & \\
X & X \otimes I \ar[l]_{s} \ar[r]^{h}  & Y \\
& X \ar[ul]^{1} \ar[ur]_{g} \ar[u]^{d_{1}} &
}
\end{equation*}
for some choice of cylinder object $X \otimes I$, in which $s$ is a
weak equivalence and $h$ is the homotopy. Then there are morphisms
\begin{equation*}
(1_{X},f) \to (s,h) \leftarrow (1_{X},g)
\end{equation*}
in $H(X,Y)$, so that $(1_{X},f)$ and $(1_{X},g)$ are in the same path
component of $H(X,Y)$. It follows that the assignment $f \mapsto
[(1_{X},f)]$ defines a function
\begin{equation*}
\psi: \pi(X,Y) \to \pi_{0}H(X,Y)
\end{equation*}
where $\pi(X,Y)$ denotes left homotopy classes of maps from $X$ to $Y$.

Theorem \ref{th 2} is a formal consequence of the following two
results:

\begin{lemma}\label{lem 3} 
Suppose that $\alpha: X \to X'$ and $\beta: Y \to Y'$
are weak equivalences. Then the induced function 
\begin{equation*}
(\alpha,\beta)_{\ast}: \pi_{0}H(X,Y) \to \pi_{0}H(X',Y')
\end{equation*}
is a bijection.
\end{lemma}

\begin{lemma}\label{lem 4} 
Suppose that $X$ is cofibrant and $Y$ is fibrant. Then
the function
\begin{equation*}
\phi: \pi_{0}H(X,Y) \to [X,Y]
\end{equation*}
is a bijection.
\end{lemma}

\begin{proof}[Proof of Lemma \ref{lem 3}]
Identify the object 
$(f,g) \in H(X',Y')$ with a map $(f,g): Z \to X'
\times Y'$ of $\mathcal{M}$ such that $f$ is a weak equivalence.

There is a factorization
\begin{equation*}
\xymatrix{
Z \ar[r]^{j} \ar[dr]_{(f,g)} & W \ar[d]^{(p_{X'},p_{Y'})} \\
& X' \times Y'
}
\end{equation*}
such that $j$ is a trivial cofibration and $(p_{X'},p_{Y'})$ is a
fibration. Observe that $p_{X'}$ is a weak equivalence. 

Form the pullback
\begin{equation*}
\xymatrix{
W_{\ast} \ar[r]^{(\alpha\times\beta)_{\ast}} \ar[d]_{(p^{\ast}_{X},p^{\ast}_{Y})} 
& W \ar[d]^{(p_{X'},p_{Y'})} \\
X \times Y \ar[r]_{\alpha \times \beta} & X' \times Y'
}
\end{equation*}
Then $(p^{\ast}_{X},p^{\ast}_{Y})$ is a fibration and $(\alpha
\times \beta)_{\ast}$ is a weak equivalence. 
The map $p^{\ast}_{X}$ is also a weak equivalence.

The assignment $(f,g) \mapsto (p^{\ast}_{X},p^{\ast}_{Y})$ defines a
function
\begin{equation*}
\pi_{0}H(X',Y') \to \pi_{0}H(X,Y)
\end{equation*}
which is inverse to $(\alpha,\beta)_{\ast}$. 
\end{proof}

\begin{proof}[Proof of Lemma \ref{lem 4}]
The canonical function $\pi(X,Y) \to [X,Y]$ is a bijection since $X$
is cofibrant and $Y$ is fibrant, and there is a commutative diagram
\begin{equation*}
\xymatrix{
\pi(X,Y) \ar[r]^-{\psi} \ar[dr]_{\cong} & \pi_{0}H(X,Y) \ar[d]^{\phi} \\
& [X,Y]
}
\end{equation*}

It suffices to show that $\psi$ is surjective, or that any object $X
\xleftarrow{f} Z \xrightarrow{g} Y$ is in the path component of some a
pair $X \xleftarrow{1} X \xrightarrow{k} Y$ for some map $k$.  Form
the diagram
\begin{equation*}
\xymatrix@R=5pt{
& Z \ar[dl]_{f} \ar[dr]^{g} \ar[dd]^{j} & \\
X && Y \\
& V \ar[ul]^{p} \ar@{.>}[ur]_{\theta} &
}
\end{equation*}
where $j$ is a trivial cofibration and $p$ is a fibration; the map $\theta$
exists because $Y$ is fibrant.
The object $X$ is cofibrant, so the trivial fibration $p$ has a
section $\sigma$, and there is a commutative diagram
\begin{equation*}
\xymatrix@R=5pt{
& X \ar[dl]_{1} \ar[dr]^{\theta\sigma} \ar[dd]^{\sigma} & \\
X && Y \\
& V \ar[ul]^{p} \ar@{.>}[ur]_{\theta} &
}
\end{equation*}
The composite $\theta\sigma$ is the required map $k$. 
\end{proof}

\section{Torsors}

\subsection{Torsors for sheaves of groups}

Suppose that $G$ is a sheaf of groups on a small Grothendieck site
$\mathcal{C}$.

A $G$-torsor is usually defined to be a sheaf $X$ with a free
$G$-action such that the map $X/G \to \ast$ is an isomorphism in the
sheaf category.  

The $G$-action on $X$ is free if and only if the the
canonical simplicial sheaf map $EG \times_{G} X \to X/G$ is a local
weak equivalence. One sees this by noting that the fundamental groups
of the Borel construction $EG \times_{G} X$ are isotropy subgroups for
the $G$-action. Further, $EG \times_{G} X$ is the nerve of a groupoid
so there are no higher homotopy groups.

It follows that a sheaf $X$ with $G$-action is a $G$-torsor if and
only if the canonical simplicial sheaf map $EG \times_{G} X \to \ast$
is a local weak equivalence.
Write $G-\mathbf{Tors}$ for the groupoid of
$G$-torsors and $G$-equivariant maps.  

Suppose given a cocycle
\begin{equation*}
\ast \xleftarrow{\simeq} Y \xrightarrow{\alpha} BG
\end{equation*}
in the category of simplicial sheaves. Form pullback
\begin{equation}\label{eq 1}
\xymatrix{
\pb(Y) \ar[r] \ar[d] & Y \ar[d]^{\alpha} \\
EG \ar[r]_{\pi} & BG
}
\end{equation}
Then the simplicial sheaf $\pb(Y)$ inherits a $G$-action from the
$G$-action on $EG$, and the induced map
$EG \times_{G} \pb(Y) \to Y$ 
is a weak equivalence. The square (\ref{eq 1}) is locally homotopy
cartesian, and it follows that the map $\pb(Y) \to
\tilde{\pi}_{0}\pb(Y)$ is a $G$-equivariant weak equivalence.  Here,
$\tilde{\pi}_{0}\pb(Y)$ is the sheaf of path components of the
simplicial sheaf $\pb(Y)$, in this case identified with a constant
simplicial sheaf.  

It follows that the maps
\begin{equation*}
EG \times_{G} \tilde{\pi}_{0}\pb(Y) \leftarrow EG \times_{G} \pb(Y) 
\to Y \simeq \ast
\end{equation*}
are weak equivalences, so that $\tilde{\pi}_{0}\pb(Y)$ is a
$G$-torsor.
A functor
\begin{equation}\label{eq 2}
H(\ast,BG) \to G-\mathbf{Tors}
\end{equation}
is therefore defined by sending the cocycle $\ast \xleftarrow{\simeq}
Y \to BG$ to the torsor $\tilde{\pi}_{0}\pb(Y)$.
 
A functor
\begin{equation}\label{eq 3}
G-\mathbf{Tors} \to H(\ast,BG)
\end{equation}
is defined by sending a $G$-torsor $X$ to the cocycle $\ast
\xleftarrow{\simeq} EG \times_{G} X \to BG$.

\begin{theorem}\label{th 5} 
Suppose that $\mathcal{C}$ is a small Grothendieck site, and that $G$
is a sheaf of groups on $\mathcal{C}$. Then the functors (\ref{eq 2})
and (\ref{eq 3}) induce bijections
\begin{equation*}
[\ast,BG] \cong \pi_{0}H(\ast,BG) \cong \pi_{0}(G-\mathbf{Tors}) =
H^{1}(\mathcal{C},G).
\end{equation*}
\end{theorem}

\begin{proof}
The functors (\ref{eq 2}) and (\ref{eq 3}) induce functions
$\pi_{0}H(\ast,BG) \to \pi_{0}(G-\mathbf{Tors})$ and
$\pi_{0}(G-\mathbf{Tors}) \to \pi_{0}H(\ast,BG)$ which are inverse to
each other. To see this, there are two statements to verify, only one
of which is non-trivial.

To show that
the composite
\begin{equation*}
\pi_{0}(G-\mathbf{Tors}) \to \pi_{0}H(\ast,G) \to \pi_{0}(G-\mathbf{Tors})
\end{equation*}
is the identity, one uses the fact that the diagram
\begin{equation*}
\xymatrix{
X \ar[r] \ar[d] & EG \times_{G} X \ar[d] \\
\ast \ar[r] & BG
}
\end{equation*} 
is locally homotopy cartesian for each $G$-sheaf $X$. This is a
consequence of Quillen's Theorem B \cite[IV.5.2]{GJ}. 
\end{proof}

The identification of the non-abelian cohomology invariant
$H^{1}(\mathcal{C},G)$ with morphisms $[\ast,BG]$ in the homotopy
category of simplicial sheaves of Theorem \ref{th 5} is, at this
writing, almost twenty years old \cite{UH}. Unlike the original proof,
the demonstration given here contains no references to hypercovers or
pro objects.

\subsection{Diagrams and torsors}

There is a local model structure for simplicial presheaves on a small site
$\mathcal{C}$ which is Quillen equivalent to the local model structure
for simplicial sheaves \cite{J2} --- this has also been known for some
time. Just recently  \cite{J7}, it has been shown that
both the torsor concept and the homotopy classification result Theorem
\ref{th 5} admit substantial generalizations, to the context of
diagrams of simplicial presheaves defined on presheaves of categories
enriched in simplicial sets.

To explain a little, when I say that $A$ is a presheaf of categories
enriched in simplicial sets I mean the standard thing, 
that $A$ consists of a presheaf $\Ob(A)$ and a simplicial
presheaf $\Mor(A)$, together with source and target maps $s,t: \Mor(A)
\to \Ob(A)$, a map $e:\Ob(A) \to \Mor(A)$ which is a section for both
$s$ and $t$, and an associative law of composition 
\begin{equation*}
\Mor(A)
\times_{t,s} \Mor(A) \to \Mor(A)
\end{equation*}
for which the map $e$ is a two-sided
identity. Here, the simplicial presheaf $\Mor(A) \times_{t,s} \Mor(A)$
is defined by the pullback
\begin{equation*}
\xymatrix{
\Mor(A) \times_{t,s} \Mor(A) \ar[r] \ar[d] & \Mor(A) \ar[d]^{s} \\
\Mor(A) \ar[r]_{t} & \Ob(A)
}
\end{equation*}
and describes composable pairs of morphisms in $A$.  To say it a
different way, $A$ is a category object in simplicial presheaves with
simplicially discrete objects.

An $A$-diagram $X$ (expressed internally \cite[p.325]{Bor},
\cite[p.240]{MM}) consists of a simplicial presheaf map $\pi: X \to
\Ob(A)$, together with an action
\begin{equation}\label{eq 4}
\xymatrix{
\Mor(A) \times_{s,\pi} X \ar[r]^-{m} \ar[d] & X \ar[d]^{\pi} \\
\Mor(A) \ar[r]_{t} & \Ob(A)
}
\end{equation}
of the morphism object $\Mor(A)$ which is associative and respects
identities.

The proof of Theorem \ref{th 5} uses Quillen's Theorem B, which can be
interpreted as saying that if $Y: I \to s\mathbf{Set}$ is an ordinary
diagram of simplicial sets defined on a small category $I$ such that
every morphism $i \to j$ induces a weak equivalence $Y_{i} \to Y_{j}$,
then the pullback diagram
\begin{equation*}
\xymatrix{
Y \ar[r] \ar[d]_{\pi} & \hocolim_{I} Y \ar[d] \\
\Ob(I) \ar[r] & BI
}
\end{equation*}
is homotopy cartesian. Here, we write $Y = \bigsqcup_{i \in \Ob(I)}
Y_{i}$. This result automatically holds when the index category $I$ is
a groupoid, but for more general index categories we have to be more
careful.

Say that $Y: I \to s\mathbf{Set}$ is a {\it diagram of equivalences} if all
induced maps $Y_{i} \to Y_{j}$ are weak equivalences of simplicial
sets, and observe that this is equivalent to the requirement that the
corresponding action
\begin{equation*}
\xymatrix{
\Mor(I)\times_{s,\pi} Y \ar[r]^-{m} \ar[d] & Y \ar[d]^{\pi} \\
\Mor(I) \ar[r]_{t} & \Ob(I)
}
\end{equation*}
is homotopy cartesian. This is the formulation that works in general:
given a presheaf of categories $A$ enriched in simplicial sets, I say
that an $A$-diagram $X$ is a diagram of equivalences if the action
diagram (\ref{eq 4}) is homotopy cartesian for the local model
structure on simplicial presheaves.

The Borel construction $EG \times_{G} Z$ for a sheaf (or presheaf) $Z$
having an action by a sheaf of groups $G$ is the homotopy colimit for
the action, thought of as a diagram defined on $G$. In general, 
say that an $A$-diagram $X$ is an {\it $A$-torsor} if
\begin{itemize}
\item[1)] $X$ is a diagram of equivalences, and
\item[2)] the canonical map $\hocolim_{A} X \to \ast$ is a local weak
equivalence.
\end{itemize}
A map $X \to Y$ of $A$-torsors is just a natural transformation,
meaning a simplicial presheaf map
\begin{equation*}
\xymatrix@C=8pt{
X \ar[rr]^{f} \ar[dr] && Y \ar[dl] \\
& \Ob(A)
}
\end{equation*}
over $\Ob(A)$ which respects actions in the obvious sense. Write
$A-\mathbf{Tors}$ for the corresponding category of $A$-torsors. This
category of torsors is not a groupoid in general, but one can show
that every map of $A$-torsors is a local weak equivalence.

Now here's the theorem:

\begin{theorem}\label{th 6}
Suppose that $A$ is a presheaf of categories enriched on simplicial
sets on a small Grothendieck site $\mathcal{C}$. Then the homotopy
colimit functor induces a bijection
\begin{equation*}
\pi_{0}(A-\mathbf{Tors}) \cong \pi_{0}H(\ast,BA) \cong [\ast,BA].
\end{equation*}
\end{theorem}
This result is proved in \cite{J7}, by a method which generalizes the
proof of Theorem \ref{th 5}. This same collection of ideas is also
strongly implicated in the homotopy invariance results for stack
cohomology which appear in \cite{J6.5}.

The definition of $A$-torsor and the homotopy classification result
Theorem \ref{th 6} have analogues in localized model categories
of simplicial presheaves, provided that those model structures are
right proper (so that Theorem \ref{th 2} applies). This is proved
in \cite{J7}. The motivic model category of Morel and Voevodsky
\cite{MV} is an example of a localized model structure for which this
result holds.

\subsection{Stack completion}

The overall technique displayed in this section specializes to give an
explicit model for the stack associated to a sheaf of groupoids $H$
--- see \cite{J8}. In general, the stack associated to $H$ has global
sections with objects given by the ``discrete'' $H$-torsors. A
discrete $H$-torsor is a $H$-torsor $X$ as above, with the extra
requirement that $X$ is simplicially constant on a sheaf of vertices;
alternatively, one could say that $X$ is a $H$-functor taking values
in sheaves such that the map $\hocolim_{H}X \to \ast$ is a weak
equivalence. This construction of the stack associated to $H$ is a
direct generalization of the classical observation that the groupoid
of $G$-torsors forms the stack associated to a sheaf of groups $G$.

\subsection{Homotopy colimits}

Suppose that $G$ is a presheaf of groupoids.

Write $s\Pre(\mathcal{C})^{G}$ for the category of $G$-diagrams in
simplicial presheaves, defined for $G=A$ as above. Following
\cite{J6.5}, this category of $G$-diagrams has an injective model
structure for which the weak equivalences (respectively cofibrations)
are those maps of $G$-diagrams
\begin{equation*}
\xymatrix@C=8pt{
X \ar[rr]^{f} \ar[dr] && Y \ar[dl] \\
& \Ob(G)
}
\end{equation*}
for which the simplicial presheaf map $f: X \to Y$ is a local weak
equivalence (respectively cofibration).

Write $s\Pre(\mathcal{C})/BG$ for the category of simplicial presheaf
maps $X \to BG$, with morphisms
\begin{equation*}
\xymatrix@C=8pt{
X \ar[rr]^{g} \ar[dr] && Y \ar[dl] \\
& BG
}
\end{equation*}
In a standard way, this category inherits a model structure from
simplicial presheaves, for which a map as above is a weak equivalence
(respectively cofibration, fibration) if and only if the simplicial
presheaf map $g: X \to Y$ is a local weak equivalence (respectively
cofibration, global fibration) of simplicial presheaves.

The homotopy colimit construction defines a functor 
\begin{equation*}
\hocolim_{G}: s\Pre(\mathcal{C})^{G} \to s\Pre(\mathcal{C})/BG.
\end{equation*}
Since $G$ is a presheaf of groupoids, this functor has a left adjoint
\begin{equation*}
\pb: s\Pre(\mathcal{C})/BG \to s\Pre(\mathcal{C})^{G},
\end{equation*}
which is defined at $X \to BG$ in sections for $U \in \mathcal{C}$ by
pulling back the map $X(U) \to BG(U)$ over the canonical maps $BG(U)/x
\to BG(U)$. The functor $\pb$ preserves cofibrations and weak
equivalences (this by Quillen's Theorem B), but more is true: the
functors $\pb$ and $\hocolim_{G}$ define a Quillen equivalence
\begin{equation*}
\pb: s\Pre(\mathcal{C})/BG \leftrightarrows s\Pre(\mathcal{C})^{G}:
\hocolim_{G}
\end{equation*} 
relating the model structures for these categories.
This result is Lemma 18 of \cite{J6.5}.

Suppose now that $X$ is a fixed $G$-diagram in simplicial presheaves,
and write $G-\mathbf{Tors}/X$ for the category with objects consisting
of $G$-diagram morphisms $A \to X$ with $A$ a $G$-torsor. The
morphisms of $G-\mathbf{Tors}/X$ are commutative diagrams
\begin{equation*}
\xymatrix@R=8pt{
A \ar[dr] \ar[dd]_{\theta} & \\
& X \\
B \ar[ur]
}
\end{equation*}
of morphisms of $G$-diagrams. One can show that the map $\theta$ must
be a weak equivalence of $G$-diagrams in all such pictures. The
homotopy colimit functor defines a functor
\begin{equation*}
\hocolim_{G}: G-\mathbf{Tors}/X \to H(\ast,\hocolim_{G} X)
\end{equation*}
in an essentially obvious way: the $G$-diagram map $A \to X$ is taken
to the cocycle
\begin{equation*}
\ast \xleftarrow{\simeq} \hocolim_{G}A \to \hocolim_{A}X.
\end{equation*}
On the other hand, given a cocycle $\ast \xleftarrow{\simeq} Y \to
\hocolim_{G}X$ the adjoint $\pb(Y) \to X$ is an object of
$G-\mathbf{Tors}/X$, so the adjunction defines a functor
\begin{equation*}
\pb: H(\ast,\hocolim_{G}X) \to G-\mathbf{Tors}/X.
\end{equation*} 
These functors are adjoint and therefore define inverse functions in
path components, so we have the following:

\begin{theorem}\label{th 7.5}
Suppose that $G$ is a presheaf of groupoids and that $X$ is a
$G$-diagram in simplicial presheaves. Then there are natural
bijections
\begin{equation*}
\pi_{0}(G-\mathbf{Tors}/X) \cong \pi_{0}H(\ast,\hocolim_{G}X) 
\cong [\ast,\hocolim_{G} X].
\end{equation*}
\end{theorem}

The identification of $\pi_{0}H(\ast,\hocolim_{G}X)$ with morphisms
$[\ast,\hocolim_{G}X]$ in the homotopy category of simplicial
presheaves in the statement of Theorem \ref{th 7.5} is a consequence
of Theorem \ref{th 2}.

Theorem \ref{th 7.5} is a generalization of Theorem 16 of \cite{J6},
which deals with the case where $G$ is a sheaf of groups and $X$ is a
sheaf (ie. constant simplicial sheaf) with $G$-action.  

\begin{example}
Let $et\vert_{S}$ be the \'etale site for a scheme $S$,
and suppose that $G$ is a discrete group acting by permutations on a
torus $\mathbb{G}_{m}^{\times n}$ defined over $S$ (more concretely,
think of the action of the Weyl group on a maximal torus in an
algebraic group). The classifying object $B(G \ltimes
\mathbb{G}_{m}^{\times n})$ is weakly equivalent to the homotopy
colimit $EG \times_{G} B(\St\mathbb{G}_{m})^{\times n}$, where
$\St\mathbb{G}_{m}$ is the stack completion for the multiplicative
group $\mathbb{G}_{m}$. The presheaf of groupoids $\St\mathbb{G}_{m}$
is equivalent to the groupoid of line bundles in each section. It
follows from Theorem \ref{th 7.5} that every torsor for the group $G \ltimes
\mathbb{G}_{m}^{\times n}$ can be identified with a choice of
family of line bundles $L_{U} = (L_{1},L_{2}, \dots ,L_{n})$, one for
each element $U \to S$ of an \'etale covering family for $S$, together
with group elements $g_{U,V} \in G$ inducing permutation isomorphisms
$L_{U}\vert_{U \cap V} \cong L_{V}\vert_{U \cap V}$,
where the group elements $g_{U,V} \in G$ form a cocycle $C(U) \to
BG$ defined on the \v{C}ech resolution associated to the covering
family $U \to S$. 
\end{example}

\section{Abelian sheaf cohomology}

Suppose that $A$ is a sheaf of abelian groups, and let $A \to J$ be an
injective resolution of $A$, thought of as a $\mathbb{Z}$-graded chain
complex and concentrated in negative degrees. Identify $A$ with a
chain complex concentrated in degree $0$, and consider the shifted
chain map $A[-n] \to J[-n]$. Observe that $A[-n]$ is the chain complex
consisting of $A$ concentrated in degree $n$. Recall that $K(A,n) =
\Gamma A[-n]$ defines the Eilenberg-Mac Lane sheaf associated to $A$,
where $\Gamma: \mathbf{Ch}_{+} \to s\mathbf{Ab}$ is the functor
appearing in the Dold-Kan correspondence. Let $K(J,n) = \Gamma
T(J[-n])$ where $T(J[-n])$ is the good truncation of $J[-n]$ in
non-negative degrees. In particular, $T(J[-n])_{0}$ is the kernel of
the boundary map $J_{-n} \to J_{-n-1}$.

As usual, write $\mathbb{Z}X$ for the free simplicial abelian group on
a simplicial set $X$, and write $NA$ for the normalized chain complex
of a simplicial abelian group $A$. Write $\pi_{ch}(C,D)$ for chain
classes of maps between presheaves of $\mathbb{Z}$-graded
chain complexes $C$ and $D$.

\begin{lemma}\label{lem 7}
Every local weak equivalence of simplicial presheaves $f: X \to Y$
induces an isomorphism
\begin{equation*}
\pi_{ch}(N\mathbb{Z}Y,J[-n]) \xrightarrow{\cong}
\pi_{ch}(N\mathbb{Z}X,J[-n])
\end{equation*}
in chain homotopy classes for all $n \geq 0$.
\end{lemma}

\begin{proof} 
Starting with the third quadrant bicomplex
$\hom(\mathbb{Z}X_{p},J^{q})$ one constructs a spectral sequence
\begin{equation*}
E_{2}^{p,q} = \Ext^{q}(\tilde{H}_{p}(X),A) \Rightarrow
\pi_{ch}(N\mathbb{Z}X,J[-p-q])
\end{equation*}
(see \cite{J1}).
The map $f$ induces a homology sheaf isomorphism
$N\mathbb{Z}X \to N\mathbb{Z}Y$, and then a comparison
of spectral sequences
gives the desired result. 
\end{proof}

Recall \cite{BK} that the category of $\mathcal{C}^{op}$-diagrams in
simplicial sets has a projective model structure for which the
fibrations (respectively weak equivalences) are the maps $f: X \to Y$
which are defined sectionwise (aka. pointwise) in the sense that each
map $f:X(U) \to Y(U)$, $U \in \Ob(\mathcal{C})$ is a fibration
(respectively weak equivalence) of simplicial sets.

If two chain maps $f,g: N\mathbb{Z}X \to J[-n]$ are chain homotopic,
then the corresponding maps $f_{\ast},g_{\ast}: X \to K(J,n)$ are
right homotopic in the projective model structure for
$\mathcal{C}^{op}$-diagrams.  Choose a sectionwise trivial fibration
$\pi: W \to X$ such that $W$ is projective cofibrant. Then
$f_{\ast}\pi$ is left homotopic to $g_{\ast}\pi$ for some choice of
cylinder object $W \otimes I$ for $W$, again in the projective
structure. This means that there is a diagram
\begin{equation*}
\xymatrix{
& W \ar[r]^{\pi} \ar[dl]_{1} \ar[d]^{i_{0}} & X \ar[dr]^{f_{\ast}} & \\
W & W \otimes I \ar[rr]^{h} \ar[l]_{s} && K(J,n) \\
& W  \ar[r]_{\pi} \ar[ul]^{1} \ar[u]_{i_{1}} & X \ar[ur]_{g_{\ast}} &
}
\end{equation*}
where the maps $s, i_{0}, i_{1}$ are all part of the cylinder object
structure for $W \otimes I$, and are sectionwise weak equivalences. It
follows that
\begin{equation*}
(1,f_{\ast}) \sim (\pi,f_{\ast}\pi) \sim (\pi s,h) \sim
(\pi,g_{\ast}\pi) \sim (1,g_{\ast})
\end{equation*}
in $\pi_{0}H(X,K(J,n))$, where $H(X,K(J,n))$ is a cocycle category for
the local model structure on simplicial presheaves.  It follows that
there is a well defined abelian group homomorphism
\begin{equation*}
\phi:\pi_{ch}(N\mathbb{Z}X,J[-n]) \to \pi_{0}H(X,K(J,n)).
\end{equation*}
This map is natural in simplicial presheaves $X$.

\begin{lemma}\label{lem 8}
The map $\phi$ is an isomorphism.
\end{lemma}

\begin{proof}
Suppose that $X \xleftarrow{f} Z \xrightarrow{g} K(J,n)$ is an object
of $H(X,K(J,n))$. The map $f$ is a local weak equivalence, so by Lemma
\ref{lem 7} there is a unique chain homotopy class $[v]: N\mathbb{Z}X
\to J[-n]$ such that $[v_{\ast}f] = [g]$.
This chain homotopy class $[v]$ is independent of
the choice of representative for the component of $(f,g)$. We
therefore have a well defined function
\begin{equation*}
\psi:\pi_{0}H(X,K(J,n)) \to \pi_{ch}(N\mathbb{Z}X,J[-n]).
\end{equation*}
The composites $\psi \cdot \phi$ and $\phi \cdot
\psi$ are identity morphisms. 
\end{proof}

\begin{corollary}\label{cor 9}
Suppose that $A$ is a sheaf of abelian groups on
$\mathcal{C}$, and let $A \to J$ be an injective resolution of
$A$ in the category of abelian sheaves. 
\begin{itemize}
\item[1)]
Let $X$ be a simplicial
presheaf on $\mathcal{C}$. Then there is a natural isomorphism
\begin{equation*}
\pi_{ch}(N\mathbb{Z}X,J[-n]) \cong [X,K(A,n)].
\end{equation*}
\item[2)]
There is a natural isomorphism
\begin{equation*}
H^{n}(\mathcal{C},A) \cong [\ast,K(A,n)].
\end{equation*}
relating sheaf cohomology to morphisms in the homotopy category of
simplicial presheaves (or sheaves).
\end{itemize}
\end{corollary}

\begin{proof}
This result is a consequence of Theorem \ref{th 2} and Lemma \ref{lem
8}. Observe that the map $K(A,n) \to K(J,n)$ is a local weak
equivalence.

The second statement is a consequence of the
first, and arises from the case where $X$ is the terminal simplicial
presheaf $\ast$.  
\end{proof}

\section{Group extensions and $2$-groupoids}

In this section we shall see that group extensions can be classified
by path components of cocycles in $2$-groupoids, by a very simple
argument. 

This is subject to knowing a few things about $2$-groupoids and their
homotopy theory.  Recall that a $2$-groupoid $H$ is a groupoid
enriched in groupoids. Equivalently, $H$ is a groupoid enriched in
simplicial sets such that all simplicial sets of morphisms $H(x,y)$
are nerves of groupoids. The object $H$ is, in particular, a
simplicial groupoid, and therefore has a bisimplicial nerve $BH$ with
associated diagonal simplicial set $dBH$.

One says that a map $G \to H$ of $2$-groupoids is a weak equivalence
if it induces a weak equivalence of simplicial sets $dBG \to
dBH$. There is a natural weak equivalence of simplicial sets
\begin{equation*}
dBH \simeq \overline{W}H 
\end{equation*}
relating $dBH$ to the space of universal cocycles $\overline{W}H$
\cite[V.7]{GJ}, \cite{JL}, so that $G \to H$ is a weak equivalence of
$2$-groupoids if and only if the induced map $\overline{W}G \to
\overline{W}H$ is a weak equivalence of simplicial sets.  One can also
show that a map $G \to H$ is a weak equivalence if and only if it
induces an isomorphism $\pi_{0}G_{0} \cong \pi_{0}H_{0}$ of path
components, and a weak equivalence of groupoids $G(x,y) \to
H(f(x),f(y))$ for all objects $x,y$ of $G$. Here, $G_{0}$ is the
groupoid of $0$-cells and $1$-cells of $G$, or equivalently the
groupoid in simplicial degree $0$ for the corresponding groupoid
enriched in simplicial sets.

The weak equivalences of $2$-groupoids are part of a general picture:
there is a model structure on groupoids enriched in simplicial sets,
due to Dwyer and Kan (see \cite[V.7.6, V.7.8]{GJ}), for which a map
$f: G \to H$ is a weak equivalence (respectively fibration) if and
only if the induced map $f_{\ast}: \overline{W}G \to \overline{W}H$ is
a weak equivalence (respectively fibration) of simplicial sets. This
model structure is Quillen equivalent to the standard model structure
for simplicial sets \cite[V.7.11]{GJ}.  The model structure for
groupoids enriched in simplicial sets restricts to a model structure
for $2$-groupoids, having the same definitions of fibration and weak
equivalence \cite{L}, and it is easy to see that both model structures
are right proper.

Here are some simple examples of $2$-groupoids:
\begin{itemize}
\item[1)]
If $K$ is a group, then
there is a $2$-groupoid $\mathbf{Aut}(K)$ with a single $0$-cell, with
$1$-cells given by the automorphisms of $K$, and with $2$-cells given by
homotopies (aka. conjugacies) of automorphisms of $K$.

\item[2)] Suppose that $p: G \to H$ is a surjective group homomorphism.
Then there is a
$2$-groupoid $\tilde{p}$ with a single $0$-cell, $1$-cells given by
the morphisms of $G$, and there is a $2$-cell $g \to h$ if and only if
$p(g) = p(h)$.
\end{itemize}

\noindent
There are canonical morphisms of $2$-groupoids
\begin{equation}\label{eq 5}
H \xleftarrow{\pi} \tilde{p} \xrightarrow{F} \mathbf{Aut}(K),
\end{equation}
The map $\pi$ is $p$ on $1$-cells, and takes a $2$-cell $g \to h$ to
the identity on $p(g) = p(h)$, whereas the map $F$ takes a $1$-cell
$g$ to conjugation by $g$, and takes the $2$-cell $g \to h$ to
conjugation by $hg^{-1} \in K$. The map $\pi: \tilde{p} \to H$ is also
a weak equivalence of $2$-groupoids, since the groupoid of $1$-cells
and $2$-cells of $\tilde{p}$ is the ``\v{C}ech groupoid'' associated
associated to the underlying surjective function $G \to H$.
In general, if $f: X \to Y$ is a surjective function,
then the associated \v{C}ech groupoid has objects given by the elements
of $X$, and a unique morphism $x \to y$ if and only if $f(x)=f(y)$.

We have therefore produced a cocycle (\ref{eq 5}) in $2$-groupoids from
a short exact sequence
\begin{equation*}
e \to K \xrightarrow{i} G \xrightarrow{p} H \to e
\end{equation*}
Write $\mathbf{Ext}(H,K)$ for the usual groupoid of all such exact
sequences. The cocycle construction is natural, and defines a
functor
\begin{equation*}
\phi: \mathbf{Ext}(H,K) \to H_{\ast}(H,\mathbf{Aut}(K))
\end{equation*}
taking values in the cocycle category $H_{\ast}(H,\mathbf{Aut}(K))$ in
pointed $2$-groupoids. All objects in the cocycle (\ref{eq 5}) have
unique $0$-cells, so the maps making up the cocycle are pointed in an
obvious way.

\begin{theorem}\label{th 10}
The functor $\phi$
induces isomorphisms
\begin{equation*}
\pi_{0}\mathbf{Ext}(H,K) \cong \pi_{0}H_{\ast}(H,\mathbf{Aut}(K)) 
\cong [BH,dB\mathbf{Aut}(K)]_{\ast},
\end{equation*}
where $[\ ,\ ]_{\ast}$ denotes morphisms in the pointed homotopy category.
\end{theorem}

\begin{proof}
If
\begin{equation*}
H \xleftarrow{\pi} A \xrightarrow{F}  \mathbf{Aut}(K)
\end{equation*}
is a pointed cocycle, then the base point $x \in A_{0}$ determines a
$2$-groupoid equivalence $A(x,x) \to H$.
The cocycle $F$ can therefore be
canonically replaced by its restriction to $A(x,x)$ at the base point
$x$, and the $2$-groupoid $A(x,x)$ can be identified with
a $2$-groupoid $p_{\ast}$ arising from a surjective group homomorphism
$p: L \to H$ with $2$-cells consisting of pairs $(g,h)$ such that
$p(g) = p(h)$.

Suppose given a cocycle 
\begin{equation*}
H \xleftarrow{\pi} p_{\ast} \xrightarrow{F} \mathbf{Aut}(K)
\end{equation*}
where $\pi: p_{\ast} \to H$ is determined a surjective group 
homomorphism $p: L \to H$ as above. There is a
group $E_{F}(p)$ which is the set of equivalence classes of pairs
$(k,x)$, $x \in L, k \in K$ such that
$(k,x) \sim (k',x')$
if and only if $p(x) = p(x')$ and $k' = F(x,x')k$. Recall that
$F(x,x')$ is a homotopy of the automorphisms $F(x), F(x')$ of $K$, and
is therefore defined by conjugation by an element $F(x,x') \in K$. The
product is defined by
\begin{equation*}
[(k,x)]\cdot [(l,y)] = [(kF(x)(l),xy)]
\end{equation*}
and there is a short exact sequence
\begin{equation*}
e \to K \to E_{F}(p) \to H \to e
\end{equation*} 
where $k \mapsto [(k,e)]$ and $[(k,x)] \mapsto p(x)$.

We have, with these constructions, described a functor
\begin{equation*}
\psi: H_{\ast}(H,\mathbf{Aut}(K)) \to \mathbf{Ext}(H,K).
\end{equation*}
One shows that the associated function $\psi_{\ast}$ on path
components is the inverse of the function
\begin{equation*}
\phi_{\ast}: \pi_{0}\mathbf{Ext}(H,K) \to \pi_{0}H_{\ast}(H,\mathbf{Aut}(K)).
\end{equation*}

The homotopy category of groupoids enriched in simplicial sets is
equivalent to the homotopy category of simplicial sets, and this
equivalence is induced by the universal cocycles functor
$\overline{W}$. It follows that there is a bijection
\begin{equation*}
[H,\mathbf{Aut}(K)]_{\ast} \cong [BH,dB\mathbf{Aut}(K)]_{\ast},
\end{equation*}
where the morphisms on the left are in the pointed homotopy category
of groupoids enriched in simplicial sets. The set
$[H,\mathbf{Aut}(K)]_{\ast}$ can be also identified with morphisms in
the pointed homotopy category of $2$-groupoids. One sees this most
effectively by observing that for every cocycle
\begin{equation*}
H \xleftarrow{\simeq} B \to \mathbf{Aut}(K)
\end{equation*}
in groupoids enriched in simplicial sets, the object $B$ is weakly
equivalent to its fundamental groupoid, and so the cocycle can be
canonically replaced by a cocycle in $2$-groupoids.
\end{proof}

\section{Classification of gerbes}

A gerbe is a stack $G$ which is locally path connected in the sense
that the sheaf of path components $\tilde{\pi}_{0}(G)$ is isomorphic
to the terminal sheaf. Stacks are really just homotopy types of
presheaves (or sheaves) of groupoids \cite{JT0}, \cite{J6}, \cite{H},
so one may as well say that a gerbe is a locally connected presheaf
of groupoids.

A morphism of gerbes is a morphism $G \to H$ of presheaves of
groupoids which is a weak equivalence in the sense that the induced
map $BG \to BH$ is a local weak equivalence of classifying simplicial
sheaves. Write $\mathbf{gerbe}$ for the category of gerbes.

Write $\mathbf{Grp}$ for the ``presheaf'' of $2$-groupoids whose
objects are sheaves of groups, $1$-cells are isomorphisms of sheaves
of groups, and whose $2$-cells are the homotopies of isomorphisms of
sheaves of groups. The object $\mathbf{Grp}$ is not really a presheaf
of $2$-groupoids because it's too big in the sense that it does not
take values in small $2$-groupoids.

Write $H(\ast,\mathbf{Grp})$ for the category of cocycles
\begin{equation*}
\ast \xleftarrow{\simeq} A \to \mathbf{Grp}
\end{equation*}
where $A$ is a presheaf of $2$-groupoids. One can sensibly discuss
such a category, even though the object $\mathbf{Grp}$ is too big to
be a presheaf. The category $H(\ast,\mathbf{Grp})$ is not small,
and its path components do not form a set. Similarly, the path
components of the category of gerbes do not form a
set. It is, nevertheless, convenient to display the relationship
between these objects in the following statement:

\begin{theorem}\label{th 11}
There is a bijection
\begin{equation*}
\pi_{0}H(\ast,\mathbf{Grp}) \cong \pi_{0}(\mathbf{gerbe}).
\end{equation*}
\end{theorem}

The proof of Theorem \ref{th 11} is a bit technical, and appears in
\cite{J9}.  It is relatively easy to say, however, how to get a
cocycle from a gerbe $G$. Write $\tilde{G}$ for the $2$-groupoid whose
$0$-cells and $1$-cells are the objects and morphisms of $G$,
respectively, and say that there is a unique $2$-cell $\alpha \to
\beta$ between any two arrows $\alpha,\beta: x \to y$.  Then the
canonical map $\tilde{G} \to \ast$ is an equivalence. There is a map
$F(G): \tilde{G} \to \mathbf{Grp}$ which associates to $x \in G(U)$
the sheaf $G(x,x)$ of automorphisms of $x$ on $\mathcal{C}/U$,
associates to $\alpha: x \to y$ the isomorphism $G(x,x) \to G(y,y)$
defined by conjugation by $\alpha$, and associates to a $2$-cell
$\alpha \to \beta$ the homotopy defined by conjugation by
$\beta\alpha^{-1}$. This cocycle construction effectively defines the
function
\begin{equation*}
\pi_{0}(\mathbf{gerbe}) \to \pi_{0}H(\ast,\mathbf{Grp}).
\end{equation*}
A generalized Grothendieck construction is used to define its
inverse --- the construction of a group from a cocycle in
the proof of Theorem \ref{th 10} is a special case.

One can go further: the gerbes with band $L \in
H^{1}(\mathcal{C},\mathbf{Out})$ are classified by cocycles in the
homotopy fibre over $L$ of a morphism of fibrant presheaves of
$2$-groupoids approximating the map
\begin{equation*}
\mathbf{Grp} \to \mathbf{Out}.
\end{equation*}
Here, $\mathbf{Out}$ is the groupoid of outer automorphisms, or
automorphisms modulo the homotopy relation. One can also use the same
techniques to classify gerbes locally equivalent to a fixed gerbe $G$.
These results are also proved in \cite{J9}.

\begin{remark}
Suppose that $E$ is a sheaf. An $E$-gerbe is a morphism of
pre\-sheaves of groupoids $G \to E$ such that the induced map
$\tilde{\pi}_{0}G \to E$ is an isomorphism of sheaves. Write
$\mathbf{gerbe}/E$ for the category of $E$-gerbes. An $E$-gerbe is
canonically a gerbe in the category of presheaves on the site
$\mathcal{C}/E$ fibred over $E$ \cite{J6.5}. It follows that Theorem
\ref{th 11} specializes to a homotopy classification statement
\begin{equation*}
\pi_{0}(\mathbf{gerbe}/E) \cong \pi_{0}H(E,\mathbf{Grp})
\end{equation*}
for $E$-gerbes.
\end{remark}


\catcode`\@=11 \renewenvironment{thebibliography}[1]{
\@xp\section\@xp*\@xp{\refname}%
\normalfont\footnotesize\labelsep
.5em\relax\renewcommand\theenumiv{\arabic{enumiv}}\let\p@enumiv\@empty
\list{\@biblabel{\theenumiv}}{\settowidth\labelwidth{\@biblabel{#1}}%
\leftmargin\labelwidth \advance\leftmargin\labelsep
\usecounter{enumiv}}%
\sloppy \clubpenalty\@M
\widowpenalty\clubpenalty \sfcode`\.=\@m }

\def\@biblabel#1{\@ifnotempty{#1}{[#1]}}
\catcode`\@=\active


\bibliographystyle{amsalpha}

\end{document}